\documentclass[12pt]{article}
\usepackage{amsmath,amssymb,amsbsy,amsfonts,amsthm,latexsym,
                           amsopn,amstext,amsxtra,euscript,amscd}
\begin{document}

\newtheorem{theorem}{Theorem}
\newtheorem{lemma}[theorem]{Lemma}
\newtheorem{claim}[theorem]{Claim}
\newtheorem{cor}[theorem]{Corollary}
\newtheorem{prop}[theorem]{Proposition}
\newtheorem{definition}{Definition}
\newtheorem{question}[theorem]{Open Question}

\def\cA{{\mathcal A}}
\def\cB{{\mathcal B}}
\def\cC{{\mathcal C}}
\def\cD{{\mathcal D}}
\def\cE{{\mathcal E}}
\def\cF{{\mathcal F}}
\def\cG{{\mathcal G}}
\def\cH{{\mathcal H}}
\def\cI{{\mathcal I}}
\def\cJ{{\mathcal J}}
\def\cK{{\mathcal K}}
\def\cL{{\mathcal L}}
\def\cM{{\mathcal M}}
\def\cN{{\mathcal N}}
\def\cO{{\mathcal O}}
\def\cP{{\mathcal P}}
\def\cQ{{\mathcal Q}}
\def\cR{{\mathcal R}}
\def\cS{{\mathcal S}}
\def\cT{{\mathcal T}}
\def\cU{{\mathcal U}}
\def\cV{{\mathcal V}}
\def\cW{{\mathcal W}}
\def\cX{{\mathcal X}}
\def\cY{{\mathcal Y}}
\def\cZ{{\mathcal Z}}

\def\A{{\mathbb A}}
\def\B{{\mathbb B}}
\def\C{{\mathbb C}}
\def\D{{\mathbb D}}
\def\E{{\mathbb E}}
\def\F{{\mathbb F}}
\def\G{{\mathbb G}}
\def\H{{\mathbb H}}
\def\I{{\mathbb I}}
\def\J{{\mathbb J}}
\def\K{{\mathbb K}}
\def\L{{\mathbb L}}
\def\M{{\mathbb M}}
\def\N{{\mathbb N}}
\def\O{{\mathbb O}}
\def\P{{\mathbb P}}
\def\Q{{\mathbb Q}}
\def\R{{\mathbb R}}
\def\S{{\mathbb S}}
\def\T{{\mathbb T}}
\def\U{{\mathbb U}}
\def\V{{\mathbb V}}
\def\W{{\mathbb W}}
\def\X{{\mathbb X}}
\def\Y{{\mathbb Y}}
\def\Z{{\mathbb Z}}

\def\E{{\mathbf E}}
\def\Fp{\F_p}
\def \Fq {{\F}_q}

\def\fP{{\mathfrak P}}
\def\fR{{\mathfrak R}}

\def \Fqs{{\F}_{q^s}}
\def \Fpn{{\F}_{p^n}}
\def \Fqn{{\F}_{q^n}}
\def\ep{{\mathbf{e}}_p}

\def\rem {\mathrm{rem\,}}

\def\rad {\mathrm{rad\,}}

\def\scr{\scriptstyle}
\def\\{\cr}
\def\({\left(}
\def\){\right)}
\def\[{\left[}
\def\]{\right]}
\def\<{\langle}
\def\>{\rangle}
\def\fl#1{\left\lfloor#1\right\rfloor}
\def\rf#1{\left\lceil#1\right\rceil}
\def\le{\leqslant}
\def\ge{\geqslant}
\def\eps{\varepsilon}
\def\mand{\qquad\mbox{and}\qquad}

\newcommand{\comm}[1]{\marginpar{%
\vskip-\baselineskip 
\raggedright\footnotesize
\itshape\hrule\smallskip#1\par\smallskip\hrule}}

\def\xxx{\vskip5pt\hrule\vskip5pt}


\title{\bf Multiplicative
Order of Gauss Periods}

\author{
      {\sc Omran Ahmadi}  \\
      {Department of Electrical and Computer Engineering}\\
      {University of Toronto} \\
      {Toronto, Ontario, M5S 3G4, Canada} \\
{\tt oahmadid@comm.utoronto.ca} \and
{\sc Igor E.~Shparlinski}\\
{Department of Computing}\\
{Macquarie University} \\
{Sydney, NSW 2109, Australia} \\
{\tt igor@ics.mq.edu.au} \and
{\sc Jos{\'e} Felipe Voloch} \\
{Department of Mathematics}\\
{University of Texas}\\
{Austin TX 78712 USA} \\
{voloch@math.utexas.edu}
}

\date{}
\pagenumbering{arabic}

\maketitle

\begin{abstract}
We obtain a lower bound on the multiplicative order of Gauss periods
which generate normal bases over finite fields. This bound improves
the previous bound of J.~von zur Gathen and I.~E.~Shparlinski.

Nous prouvons une borne inf\'{e}rieure pour l'ordre multiplicatif des
p\'{e}riodes de Gauss g\'{e}n\'{e}rant les bases normales sur les corps finis.
Cette borne am\'{e}liore une borne ant\'{e}rieure du \`{a} J.~von zur Gathen
et I.~E.~Shparlinski.
\end{abstract}

\section{Introduction}

For a prime power $q$ we use $\F_q$ to denote the
finite field with $q$ elements.

Normal bases are a very useful notion in the theory of finite fields,,
see~\cite{BGMMVY,LN,Shp} for the definition, basic properties and references.
One of the most interesting constructions of
normal bases come from {\it Gauss periods\/},
see~\cite{FvzGS,GvzGP,vzGN1,vzGN2,vzGS1,vzGS2,vzGS3} and
references therein.
In particular, Gauss periods of type $(n,2)$ are
of special interest, which can be defined as
follows.

Let $r=2n+1$ be a prime number coprime with
$q$ and
$\beta \in \F_{q^{2n}}$ be a  primitive
$r$th root  of unity.
Then the element
\begin{equation}
\label{eq:GP}
\alpha = \beta + \beta^{-1} \in \F_{q^n}
\end{equation}
is called a Gauss period of type $(n,2)$.
The Gauss period of type $(n,2)$ can be defined for composite $r$ too,
see~\cite{FvzGS}, however we do not consider them in this
paper (neither we study Gauss period of type $(n,k)$ for $k\neq 2$).

It is well-known that the minimal polynomial of $\beta$ over $\Fq$
is of degree $t$, where $t$ is the multiplicative order of $q$ modulo
$r$. Thus $t | 2n$.

It is also well known that $\alpha$ given by~\eqref{eq:GP},
generates a normal basis of $\Fqn$ if an only if $\gcd(2n/t,n) =1$,
which, therefore, is possible if and only if
\begin{itemize}
\item  $t = 2n = r-1$, that is, $q$ is a primitive root modulo $r$;
\item  $t = n = (r-1)/2$ and $n$ is odd, that is, $q$ generates the
subgroup of quadratic residues modulo $r \equiv 3 \pmod 4$
\end{itemize}
In one direction this
follows from~\cite[Lemma~5.4 and Theorem~5.5]{BGMMVY} and in
the other
direction it follows by examining the proof of these results
see also~\cite{AHM} and~\cite{ABV}.

It is shown~\cite{vzGS1} that in the first case, that is, for $t = r-1$,
$\alpha$ is of multiplicative order
\begin{equation}
\label{eq:Old Bound}
L_n \ge 2^{\sqrt{2n} + O(1)},
\end{equation}
see also~\cite{vzGS2}. This gives an explicit example of finite
field elements of exponentially large order. Here  we use
some new arguments to improve the
bound~\eqref{eq:Old Bound}.

Recent results
of Q.~Cheng~\cite{Cheng} give polynomial
time constructions of elements
of large order for certain values of $(q,n)$.
Our construction seems to apply to different sets of pairs $(q,n)$
and complement the results of ~\cite{Cheng}. Furthermore it is
interesting to establish
tighter bounds on  the size of the multiplicative order of such
classical objects as  Gauss periods of type $(n,2)$,
especially of those which generate normal bases.

Let $P(s,v)$ be the number of integer partitions
of an integer $s$ where each part appears no more
than $v$ times, that is, the number of
solutions to the equation
$$
\sum_{j=1}^{s}u_j j = s
$$
in non-negative integers $u_1, \ldots, u_s \le v$.

\begin{theorem}
\label{thm:p-bound}
Let $p$ be the characteristic of $\F_q$
and let  $q$ be a primitive root modulo a
prime
$r = 2n+1$. Then the multiplicative order $L_n$ of $\alpha$, given
by~\eqref{eq:GP}, satisfies the bound
$$
L_n \ge P(n-1, p-1).
$$
\end{theorem}

Now we can use some standard estimates to derive
an asymptotic lower bound on $L_n$.
%

\begin{cor}
\label{cor:H-R bound}
Let $p$ be the characteristic of $\F_q$
and let  $q$ be a primitive root modulo a
prime
$r = 2n+1$.
Then, uniformly over $q$, the multiplicative order $L_n$ of $\alpha$, given
by~\eqref{eq:GP}, satisfies the bound
$$
L_n \ge \exp\(\(\pi \sqrt{\frac{2(p-1)}{3p}} + o(1)\)\sqrt{n}\),
$$
   as $n \to \infty$.
\end{cor}

Note that in the worst case (when $p =2$)
$ \exp\( \pi \sqrt{2/6}\) = 6.1337\ldots$
while $ \exp\( \pi \sqrt{2/3}\) = 13.0019\ldots$ (which
corresponds to $p \to \infty$).
On the other hand, we have $2^{\sqrt{2}} = 2.6651\ldots$.

\section{Proof of Theorem~\ref{thm:p-bound}}

Let us consider the set
$$
\fP=
\left\{(u_1,\ldots,u_{n-1})\in \Z_{\ge 0}^{n-1}\ |\
\sum_{j=1}^{n-1}u_j j = n-1, \ u_1, \ldots, u_n \le p-1
\right\}.
$$

Now, for  $j=1,2,\dots,n-1$ we define an integer $z_j$
by  $q^{z_j}\equiv j \pmod{r}$, $0 \le z_j < r$
(which is possible since $q$ is a primitive root modulo $r$).

For every partition $\cU = (u_1, \ldots, u_{n-1}) \in \fP$ we
put
$$
Q_\cU = \sum_{j=1}^{n-1} u_jq^{z_j}.
$$
We now consider the powers
$$
\alpha^{Q_\cU}  =  \prod_{j=1}^{n-1} \alpha^{u_jq^{z_j}}
= \prod_{j=1}^{n-1}\(\beta + \beta^{-1}\)^{u_jq^{z_j}} =
\prod_{j=1}^{n-1} \(\beta^{q^{z_j}} + \beta^{-q^{z_j}}\)^{u_j}
$$
taken for all  $ \cU \in \fP$.
Since $\beta^r =1$, we have
\begin{equation}
\label{eq:alpha and beta}
\alpha^{Q_\cU} =  \prod_{j=1}^{n-1} \(\beta^j + \beta^{-j}\)^{u_j} =
\beta^{-(n-1)}\prod_{j=1}^{n-1}\(\beta^{2j} +1\)^{u_j}.
\end{equation}

Clearly it suffices to show that for two  distinct partitions
$ \cU,\cV \in \fP$ we have $\alpha^{Q_\cU} \ne\alpha^{Q_\cV}$.

We now assume that there are two distinct partitions
$$
\cU= (u_1, \ldots, u_{n-1}),\ \cV= (v_1, \ldots, v_{n-1})
\in
\fP
$$ with
\begin{equation}
\label{eq:Assumption}
\alpha^{Q_\cU}=\alpha^{Q_\cV}.
\end{equation}
By~\eqref{eq:alpha and beta} we conclude that
$$
\prod_{j=1}^{n-1}\(\beta^{2j} +1\)^{u_j}=
\prod_{j=1}^{n-1}\(\beta^{2j} +1\)^{v_j}.
$$
Since the characteristic polynomial of $\beta$ is the $r$-th cyclotomic
polynomial $\varPhi_r(X)$, we obtain polynomial  divisibility
\begin{equation}
\label{eq:div ABC}
\varPhi_r(X) \mid U(X)-V(X)
\end{equation}
where
$$
U(X) = \prod_{j=1}^{n-1}\(X^{2j} +1\)^{u_j}, \qquad
V(X) = \prod_{j=1}^{n-1}\(X^{2j} +1\)^{v_j},
$$
are polynomials of degree $2(n-1) < 2n = r-1  = \deg \varPhi_r(X)$.(notice
that $r$ is a prime number and $q$ is a primitive root modulo $r$)
Hence~\eqref{eq:div ABC} implies that $U(X)= V(X)$.
After removing common factors, the identity
$$
\prod_{j=1}^{n-1}\(X^{2j} +1\)^{u_j}= \prod_{j=1}^{n-1}\(X^{2j} +1\)^{v_j}
$$
leads to the relation
\begin{equation}
\label{eq:Ident}
\prod_{h\in \cH} \(X^{2h} +1\)^{y_h}= \prod_{k \in \cK}\(X^{2k} +1\)^{z_k}
\end{equation}
for two disjoint sets $\cH,\cK \in \{1, \ldots, n-1\}$
and some positive integers $y_h$, $h\in \cH$,
and $z_k$, $k \in\cK$.
Since it is now clear that
$$
\gcd\(\prod_{h\in \cH}y_h  \prod_{k \in \cK}z_k, p\) =1,
$$
   the term $X^{2f}$ where
$f$ is the smallest element of $\cH \cup \cK$ occurs
only on one side of~\eqref{eq:Ident}, which makes this
identity impossible.

Therefore~\eqref{eq:Assumption} cannot hold
and the result follows.

\section{Proof of Corollary~\ref{cor:H-R bound}}

Unfortunately, a uniform lower bound with respect to $v$
on $P(s,v)$ does not seem to be in the literature.
However, by~\cite[Corollary~1.3]{And} we have
$$
P(s,v) = Q(s,v+1)
$$
where $Q(s,d)$ is the number of integer partitions
of an integer $s$ where each part is not divisible
by $d$, that is, the number of
solutions to the equation
$$
\sum_{j=1}^{s}u_j j = s
$$
in non-negative integers $u_1, \ldots, u_s$ such that $u_j =0$
for $j \equiv 0 \pmod d$, $j =1,\ldots, n$.

By~\cite[Corollary~7.2]{Hagis}, applied to a set
$\{1, \ldots, (\ell-1)/2\}$ for a fixed prime $\ell$
(thus $r = (\ell-1)/2$)
implies that
\begin{equation}
\label{eq:Q-bound}
Q(s,\ell) \ge  \exp\(\(\pi \sqrt{\frac{2(\ell-1)}{3\ell}} +
o(1)\)\sqrt{n}\).
\end{equation}
Therefore there is a function $\lambda(s)\to \infty$ as
$s \to \infty$, such that~\eqref{eq:Q-bound} holds
uniformly over all primes $\ell \le \lambda(s)$.

Now taking $\ell$ as the largest prime with
$$
\ell \le \min\{p, \lambda(n-1)\}
$$
we obtain
$$
P(n-1,p-1) \ge P(n-1,\ell-1) = Q(n-1,\ell).
$$
Applying~\eqref{eq:Q-bound} we obtain the desired estimate.
Indeed, if $\ell = p$ this is obvious. If
$\ell \le  \lambda(n-1) < p$ then by the prime number
theorem $\ell \sim \lambda(n-1)$ as $n \to \infty$.
Therefore,
$$
   \frac{\ell-1}{\ell} = 1 + O(1/\lambda(n-1)) \mand
   \frac{p-1}{p} = 1 + O(1/\lambda(n-1)) .
$$

\section{Remarks}

It seems to be natural to use the approach
of~\cite{Bern,Vol}, based on the polynomial
$ABC$-theorem, see~\cite{Sny}, in order to
obtain good bounds on $L_n$. In fact this
is possible indeed, however it seems to lead to a result
which is slightly weaker than the bound of
Theorem~\ref{thm:p-bound}. In fact, instead
of the set $\fP$ one seems to need to consider
sets of the shape
$$
\fR_s(N)=
\left\{(u_1,\ldots,u_{s})\in \Z_{\ge 0}^{n-1}\ |\ \sum_{j=1}^{s}u_j j = N
\right\}
$$
with $s \sim \alpha n^{1/2}$ and $N = \beta n$,
where $\alpha$ and $\beta$ are positive constants
(which are to be optimised).
We remark that an asymptotic formula
for  $\# \fR_s(N)$ is given
by a result of G.~Szekeres~\cite{Szek1,Szek2}. Using this approach we
have been able to get a stronger
bound than~\eqref{eq:Old Bound} but marginally weaker
than that of Theorem~\ref{thm:p-bound}.  Still, it seems
quite plausible that a use of the polynomial
$ABC$-theorem may lead to stronger bounds. We pose this
as an open question.

\section*{Acknowledgements}

The authors are very grateful to George Andrews for
providing several crucial references.

    This work was initiated
during very pleasant visits by I.~S.\ to
Department of Combinatorics \& Optimization of
the University of Waterloo
the Department of Mathematics of the University of Texas; the
hospitality, support and stimulating research
atmosphere of these institutions are gratefully
acknowledged. During the preparation of this paper, I.~S.\ was
supported in part by ARC grant DP0556431.

\end{document}